\pdfoutput=1
\RequirePackage{ifpdf}
\ifpdf 
\documentclass[pdftex]{sigma}
\else
\documentclass{sigma}
\fi

\usepackage[all]{xy}

\newcommand{\G}[1]{\mathfrak{#1}}
\newcommand{\C}[1]{\mathcal{#1}}
\newcommand{\B}[1]{\mathbb{#1}}
\newcommand{\lmod}[1]{{#1}\text{-{\bf Mod}}}

\newcommand{\CH}{{\rm CH}}

\newcommand{\tr}{\triangleright}
\newcommand{\tl}{\triangleleft}

\renewcommand{\geq}{\geqslant}

\numberwithin{equation}{section}

\newtheorem{Theorem}{Theorem}[section]
\newtheorem*{Theorem*}{Theorem}
\newtheorem{Corollary}[Theorem]{Corollary}
\newtheorem{Lemma}[Theorem]{Lemma}
\newtheorem{Proposition}[Theorem]{Proposition}
 { \theoremstyle{definition}

 }

\begin{document}
\allowdisplaybreaks

\newcommand{\arXivNumber}{2009.14801}

\renewcommand{\PaperNumber}{063}

\FirstPageHeading

\ShortArticleName{Birational Equivalences and Generalized Weyl Algebras}

\ArticleName{Birational Equivalences\\ and Generalized Weyl Algebras}

\Author{Atabey KAYGUN}
\AuthorNameForHeading{A.~Kaygun}

\Address{Istanbul Technical University, Istanbul, Turkey}
\Email{\href{mailto:kaygun@itu.edu.tr}{kaygun@itu.edu.tr}}
\URLaddress{\url{http://web.itu.edu.tr/~kaygun/}}

\ArticleDates{Received June 24, 2024, in final form July 21, 2025; Published online July 30, 2025}

\Abstract{We calculate suitably localized Hochschild homologies of various quantum groups and Podle\'s spheres after realizing them as generalized Weyl algebras (GWAs). We use the fact that every GWA is birationally equivalent to a smash product with a 1-torus. We also address and solve the birational equivalence problem, and the birational smoothness problem for GWAs.}


\Classification{17B37; 16E40; 14E99}

\section{Introduction}

A birational equivalence is an algebra morphism that becomes an isomorphism after a suitable
localization. In this paper, we show that every generalized Weyl algebra (GWA) is birationally
equivalent to a smash product with a rank-1 torus. This fact significantly simplifies their
representation theory, and structure problems such as the isomorphism
problem~\cite{BavulaJordan:IsoProblem,Hodges:NCDeformations, RichardSolotar:IsoProblem,
 AlvarezVivaz:IsomorphismsOfQGWAs, Tang:IsoProblem} and the smoothness
problem~\cite{Bavula:GlobalDim, Hodges:NCDeformations,
 Liu:HomologicalSmoothnessI,Liu:HomologicalSmoothnessII, SolotarAlvarezVivas:QuantumGeneralizedWeyl}, provided one replaces isomorphisms with
suitable noncommutative birational equivalences. We address and solve a relative version of the
birational equivalence problem in Section~\ref{subsect:IsoProblem}, and the birational smoothness
problem in Section~\ref{sect:separable}. We then calculate the Hochschild homology of suitably
localized examples of GWAs in Section~\ref{sect:Calculations}.\looseness=1

Generalized Weyl algebras are defined by Bavula~\cite{Bavula:GWAI, Bavula:GlobalDim},
Hodges~\cite{Hodges:NCDeformations} and Rosenberg~\cite{Rosenber:NCAlgebraicGeometry} independently
under different disguises. Their representation theory resembles that of Lie
algebras~\cite{DGO:WeightModules, Meinel:DufloTheorem} (see Section~\ref{sect:highest-weights}),
their homologies are extensively studied~\cite{FarinatiSolotarAlvarez:GeneralizedWeyAlgebras, Liu:HomologicalSmoothnessI,
 Liu:HomologicalSmoothnessII,
 SolotarAlvarezVivas:QuantumGeneralizedWeyl}, and they found diverse uses in areas such as noncommutative
resolutions of Kleinian
singularities~\cite{Boddington:ResolutionOfTypeD,crawley-boevey1998,Levy:2009} and noncommutative
geometry of various quantum spheres and lens spaces~\cite{Brzezinski:GWAs}. Apart from
noncommutative resolutions of Kleinian singularities, the class is known to contain the ordinary
rank-1 Weyl algebra $A_1$, the enveloping algebra $U(\G{sl}_2)$ and its primitive quotients, the
quantum enveloping algebra $U_q(\G{sl}_2)$, the quantum monoid $\C{O}_q(M_2)$, the quantum groups
$\C{O}_q({\rm GL}_2)$, $\C{O}_q({\rm SL}_2)$ and $\C{O}_q({\rm SU}_2)$. We reverify that the standard Podle\'s spheres
$\C{O}_q\bigl(S^2\bigr)$~\cite{Podl87} are GWAs~\cite{Liu:HomologicalSmoothnessI}, and then we show that
parametric Podle\'s spheres $\C{O}_{q,c}\bigl(S^2\bigr)$ of Hadfield~\cite{Hadfield:Podles} are also GWAs. We
finish the paper by calculating localized Hochschild homology of all of these examples.

The Hochschild homology of quantum groups $\C{O}_q({\rm GL}_n)$ and $\C{O}_q({\rm SL}_n)$ with coefficients in a~1-dimensional character coming from a modular pair in involution is calculated for every $n\geq 1$
in~\cite{KaygunSutlu:HomologyOfQuantumLinearGroups}, and with coefficients in themselves in specific
cases in~\cite{HadfKrah05, HadfKrae06, MasuNakaWata90, Ross90}. The Hochschild cohomology of the
Podle\'s sphere was studied by Hadfield~\cite{Hadfield:Podles}, and then in the context of van den
Bergh duality~\cite{vdBer98, VanDenBergh:HochschildPoicareDualityErratum} by Kr\"ahmer~\cite{Krah08}.
Both Hadfield and Kr\"ahmer use twisted Hochschild (co)homology by the Nakayama automorphism with
coefficients in themselves. In~this paper, we~only calculate the ordinary Hochschild homology of
these algebras with coefficients in themselves since one can always move to and from the ordinary
Hochschild homology and the twisted homology via isomorphisms coming from suitable cup and cap
products~\cite{GoodKrah14,Krah08}.

In this paper, we focus on GWAs, i.e., algebras that are birationally equivalent to smash products
with rank-1 tori. The higher rank generalized Weyl algebras that (conjecturally) recover enveloping
algebras of higher rank Lie algebras and their quantizations are called \emph{twisted generalized
 Weyl algebras} (TGWAs)~\cite{Hartwig:LocallyFiniteSimpleWeights,
 Hartwig:TGWAs, MPT03,MazorchukTurowska:TGWAs}. We conjecture that TGWAs are birationally equivalent to smash products with
higher rank tori, but we leave this investigation for a future paper.\looseness=-1

The celebrated Gelfand--Kirilov conjecture, on the other hand, states that the universal enveloping
algebra $U(\G{g})$ of a finite dimensional Lie algebra is birationally equivalent to a sufficiently
high rank Weyl algebra~\cite{GelfandKirilov:Conjecture}. One of the equivalent forms of the
conjecture is that $U(\G{g})$ is birationally equivalent to the smash product of a polynomial
algebra with a torus. The conjecture is known to be false in general~\cite{AOV:CounterExamplesToGK,
 CKPR11}, but is true for a large class of Lie algebras~\cite{GelfandKirilov:Conjecture,
 Joseph:GKConjecture, Hai:GKConjecture}. The quantum analogue of the conjecture
(see~\cite[pp.~19--21 and Section~II.10.4]{BrownGoodearl:QuantumGroups} and references therein) is also
known to be true many instances~\cite{AlevDumas:QuantumGK, FutornyHartwig:QuantumGK}. In the light
of our conjecture above, we believe that the universal enveloping algebra $U(\G{g})$ of a rank-$n$
semi-simple Lie algebra is birationally equivalent to the smash product of a polynomial algebra with
an $n$-torus. We also believe that the same is true for the quantum enveloping algebras
$U_q(\G{g})$ and the quantum groups~$\C{O}_q(G)$ where one replaces the $n$-torus with a quantum
$n$-torus.

{\bf Plan of the article.}
In Section~\ref{sect:prelims}, we recall some basic facts on localizations, relative homology
of algebra extensions, smash products and biproducts. In Section~\ref{sect:GWAs}, we prove
two fundamental structure theorems for GWAs in~Sections~\ref{sect:structureTheorem}
and~\ref{sect:localizationOfSmashProducts}. Then we state and solve birational equivalence
problem for GWAs in Section~\ref{subsect:IsoProblem}. In Section~\ref{sect:Homology}, we
investigate the interactions between homology, smash biproducts and noncommutative
localizations, and in Sections~\ref{sect:separable} and~\ref{sect:localizationGWA}, we state
and solve the birational smoothness problem for GWAs. Finally, we use our machinery to
calculate suitably localized Hochschild homologies of various GWAs in
Section~\ref{sect:Calculations}.

{\bf Notations and conventions.}
We fix an algebraically closed ground field $\Bbbk$ of characteristic~0, and we set the binomial
coefficients $\binom{n}{m}=0$ whenever $m>n$ or $m<0$. All unadorned tensor products $\otimes$ are
taken over $\Bbbk$. All algebras are assumed to be unital and associative, but not necessarily
commutative or finite dimensional. We use the notation $\Bbbk[X]$ for the free unital commutative
algebra generated by a set $X$, while we use $\Bbbk\{X\}$ for the free unital algebra generated by
the same set $X$. Throughout the paper, we use $\B{T}$ to denote the algebra of Laurent polynomials
$\Bbbk\big[x,x^{-1}\big]$.

\section{Preliminaries}\label{sect:prelims}

\subsection{Noncommutative localizations}

Our main reference for noncommutative localizations is \cite[Section~10]{Lam:Lectures}.

A multiplicative submonoid $S\subseteq A$ is called a \emph{right Ore set} if for every
$s\in S$ and $u\in A$
\begin{enumerate}\itemsep=0pt
\item[(i)] there are $s'\in S$ and $u'\in A$ such that $su'=us'$, and
\item[(ii)] if $su=0$ then there is $s'\in S$ such that $us'=0$.
\end{enumerate}

If $S\subseteq A$ is a right Ore set, then there is a \emph{universal} algebra $S^{-1}A$ and a
morphism of algebras~${\iota_S\colon A\to S^{-1}}A$ such that $\varphi(S)\subseteq S^{-1}A^\times$. The
morphism $\iota_S$ is universal among such $S$ inverting morphisms, where if $\varphi\colon A\to B$
satisfies $\varphi(S)\subset B^\times$, then there is a unique morphism of algebras
$\varphi'\colon S^{-1}A\to B$ with $\varphi = \varphi'\circ \iota_S$.

In the sequel, we are going to drop the requirement that $S$ is a multiplicative submonoid and
consider the conditions above within the submonoid generated by $S$. In such cases, we are still
going to use the notation $S^{-1}A$ for the localization.

\subsection{Birational equivalences}

Given two irreducible algebraic varieties $X$ and $Y$, an algebraic map $f\colon X\to Y$ is
called a~birational equivalence if $f$ restricted to an open subvariety is an
isomorphism~\cite[Section~4.2]{GriffithsHarris:AlgGeo}. For the purposes of this article, an
algebra $A$ is \emph{an irreducible noncommutative affine variety} if $A$ has no primitive
idempotents other than 0 and 1. Thus we call a morphism of algebras $\varphi\colon A\to A'$
\emph{a birational equivalence} if there are two Ore sets $S\subset A$ and $S'\subset A'$
such that $\varphi(S)\subseteq S'$ and the extension of $\varphi$ to the localization
$S^{-1}\varphi\colon S^{-1}A\to (S')^{-1}A'$ is an isomorphism of unital associative
algebras.

\subsection{Smash biproducts}

Assume $A$ and $B$ are two unital associative algebras. A $\Bbbk$-linear map $R\colon
B\otimes A\to A\otimes B$ is called a distributive law if the following diagrams of
algebras commute:
\[
 \xymatrix{
 B\otimes B\otimes A \ar[r]^{B\otimes R}\ar[d]_{\mu_B\otimes A}
 & B\otimes A\otimes B \ar[r]^{ R \otimes B}
 & A\otimes B\otimes B \ar[d]^{A\otimes \mu_B}\\
 B\otimes A \ar[rr]^R
 & & A\otimes B\\
 B\otimes A\otimes A \ar[r]_{ R \otimes B} \ar[u]^{B\otimes\mu_A}
 & A\otimes B\otimes A \ar[r]_{A\otimes R }
 & A\otimes A\otimes B, \ar[u]_{\mu_A\otimes B}
 }\qquad \xymatrix{
 & B\ar[dr]^{1\otimes B}\ar[dl]_{B\otimes 1}\\
 B\otimes A \ar[rr]^{ R }
 & & A\otimes B.\\
 & A \ar[ul]^{1\otimes A} \ar[ur]_{A\otimes 1}
 }
\]
For notational convenience, we write
\smash{$R(b\otimes a) = R_{(1)}(a)\otimes R_{(2)}(b) $} for every $a\in A$ and $b\in B$.

For a distributive law $R\colon B\otimes A\to A\otimes B$, there is a corresponding smash
biproduct algebra~$A\#_R B$ which is $A\otimes B$ as vector spaces with the multiplication
$ (a\otimes b)(a'\otimes b') = aR_{(1)}(a')\otimes R_{(2)}(b)b'$ for every $a,a'\in A$ and
$b,b'\in B$. We will use $A\# B$ instead of $A\#_R B$ if the distributive law is clear from the
context.

\subsection{Hochschild homology}

Let $A$ be a unital associative algebra, and let $M$ be an $A$-bimodule. Consider the
graded $\Bbbk$-vector space
$\CH_*(A,M) = \bigoplus_{n\geq 0} M\otimes A^{\otimes n} $ together with linear maps
$b_n\colon \CH_n(A,M)\to\CH_{n-1}(A,M)$ defined for $n\geq 1$ via
\begin{align*}
 b_n(m\otimes a_1\otimes a_n)
 = {}& m a_1\otimes a_2\otimes\cdots\otimes a_n + \sum_{i=1}^{n-1}(-1)^i m\otimes\cdots\otimes a_ia_{i+1}\otimes\cdots\otimes a_n \\
 & + (-1)^n a_n m \otimes a_1\otimes\cdots\otimes a_{n-1}.
\end{align*}
These maps satisfy $b_nb_{n+1} = 0$ for every $n\geq 1$, and we define
$H_*(A,M) = \ker(b_n)/im(b_{n+1})$. We use the notation $HH_*(A)$ for $H_*(A,A)$.

\subsection{Smooth algebras}\label{sect:smooth}

An algebra $A$ is said to be \emph{smooth} if it has finite Hochschild homological
dimension, i.e., when
\[ \operatorname{ hh.dim}(A) := \sup\{n\in\B{N}\mid H_n(A,M)\neq 0,\,M\in\lmod{A^e}\} \] is finite. We call
an algebra $A$ \emph{$m$-smooth} if $hh.dim(A)=m$, for $m\in\B{N}$.

The simplest $0$-smooth algebras are groups ring $\Bbbk[G]$ over finite groups where $|G|$
does not divide the characteristic of $\Bbbk$, and quotients of polynomial algebras
$\Bbbk[x]/\langle f(x)\rangle$ where $f(x)$ is a separable polynomial. For $m\geq 1$, the
simplest examples of $m$-smooth algebras one can~consider are the polynomial algebras
$\Bbbk[t_i\mid i=0,\ldots,m]$ and the Laurent polynomial algebras~${\Bbbk\big[t_i,t_i^{-1}\mid i=0,\ldots,m\big]}$ with $m\geq 0$, and their smash biproducts with
$0$-smooth algebras.

\subsection{Homology of smash biproducts with smooth algebras}
Now, let us recall the following result
from~\cite[Proposition~1.5]{KaygunSutlu:HomologyOfQuantumLinearGroups}.

\begin{Proposition}\label{base2}
 Assume $P$ and $Q$ are two unital algebras together with a left flat algebra morphism
 $\varphi \colon Q \to P$. Let $M$ be a $P$-bimodule. Then there is a spectral sequence
 whose first page is given by
 \[ E^1_{i,j}
 = H_j(Q, \underbrace{M\otimes_Q P\otimes_Q \cdots\otimes_Q P}_{\text{$i$-times}})
 \] that converges to the Hochschild homology $H_*(P, M)$.
\end{Proposition}
Let us write $R_{(1)}(b\otimes\mathbf{a})\otimes R_{(2)}(b\otimes\mathbf{a})$ for
$R(b\otimes \mathbf{a})$, and defined
\[ \CH_n(A,M)_B := \frac{\CH_n(A,M)}{\{\mathbf{a}\otimes m\triangleleft b -
 R_{(1)}(b\otimes\mathbf{a})\otimes R_{(2)}(b\otimes\mathbf{a})\triangleright m \mid b\in B,\,\mathbf{a}\otimes m\in
 \CH_n(A,M)\}}.
\]
\begin{Corollary}\label{base3}
 Let $A$ and $B$ be two algebras, and let $R\colon B \otimes A \to A\otimes B$ be any
 distributive law. Then for any $A\# B$-bimodule $M$ and for all $n \geq 0$ we have
$
 H_n (A\# B, M) \cong H_n(\CH_*(A, M)_B)$,
 when $B$ is $0$-smooth.
\end{Corollary}

\begin{proof}
 We set $P=A\# B$ and $Q=B$ together with $\varphi(b) = 1\otimes b$, and then we use
 Proposition~\ref{base2}.
\end{proof}

\section{Generalized Weyl algebras}\label{sect:GWAs}

\subsection{Generalized Weyl algebras}

Assume $A$ is a unital associative algebra, let $a\in Z(A)$ and $\sigma\in \operatorname{Aut}(A)$ be fixed.
Define a~new algebra $W_{a,\sigma}$ as the quotient of the free algebra generated by $A$ and
two non-commuting indeterminates $x$ and $y$ subject to the following relations
\[
 yx-a,\qquad xy-\sigma(a),\qquad xu-\sigma(u)x,\qquad y\sigma(u)-uy
\]
 for every $u\in A$. The algebra $W_{a,\sigma}$ is called \emph{generalized
 Weyl algebra}~\cite{Bavula:GWAI, BavulaJordan:IsoProblem}.

\subsection{Algebras with automorphisms}

One standard source of distributive laws is algebras with a fixed algebra automorphism or
endomorphisms. From now on, we will assume $A$ is a domain with a fixed algebra automorphism~${\sigma\in \operatorname{Aut}(A)}$. Let $\B{T}=\Bbbk[\B{Z}] = \Bbbk\big[x,x^{-1}\big]$ be the group ring of the
free abelian group on a~single generator $\B{Z}$. Now consider the smash biproduct
$B := A\# \B{T}$ coming from the distributive law $R\colon \B{T}\otimes A\to A\otimes \B{T}$
defined as
\begin{equation}\label{monoid-dist-law} R(x^n\otimes u) = \sigma^n(u)\otimes x^n
\end{equation} for every monomial $x^n\in \B{T}$ with $n\in\B{Z}$ and $u\in A$. Then $R$
defines an invertible distributive law. In order to simplify the notation, we are going
to write $ux^i$ for every monomial $u\otimes x^i$ in~$A\# \B{T}$. If there are more
than one automorphisms in the context, we are going to write $A\#_{R(\sigma)} \B{T}$
instead of $A\# \B{R}$ to emphasize which automorphism we are using.

\subsection{The ordinary Weyl algebra}\label{sect:OrdinaryWeyl}

If we consider the subalgebra of operators generated by $T(f(t)) = t f(t)$ and
$D(f(t))= f'(t)$ on the space of differentiable functions on $\B{R}$, we see that we have
the relation $[D,T] = {DT - TD = 1}$. This algebra is the motivating example of the ordinary rank-1 Weyl
algebra~$A_1$, which is defined as the $\Bbbk$-algebra defined on two non-commuting
indeterminates $x$ and $y$ subject to the relations~${xy - yx = 1}$. One can define~$A_1$ as a GWA if we let $A=\Bbbk[t]$ where we set the
distinguished element $a=t$. We define $\sigma$ to be the algebra automorphism of~$A$ given
by~${f(t) = f(t-1)}$ for every $f(t)\in A$. Then the GWA $W_{a,\sigma}$ is the ordinary Weyl
algebra~$A_1$. See~\cite[Example~2.3]{BavulaJordan:IsoProblem}.

\subsection{A structure theorem for GWAs}\label{sect:structureTheorem}

Given an arbitrary GWA $W_{a,\sigma}$, one can also realize it as a unital subalgebra of the smash
product~$A\# \B{T}$ where $\B{T}:= \Bbbk\big[x,x^{-1}\big]$ and $R\colon \B{T}\otimes A\to A\otimes \B{T}$
is defined in equation~\eqref{monoid-dist-law}. For this, we consider the monomorphism of
$\Bbbk$-algebras $\varphi\colon W_{a,\sigma}\to A\# \B{T}$ given by
$
\varphi(u) = u$, $ \varphi(x) = x$, $ \varphi(y) = ax^{-1}
$ for every $u\in A$.

\begin{Theorem}\label{thm:UnitCase}
 For every $a\in Z(A)$, the algebra $W_{a,\sigma}$ is isomorphic to the unital subalgebra
 of the smash biproduct $A\# \B{T}$ generated by $A$, $x$ and $ax^{-1}$. Hence
 $W_{a,\sigma}$ is isomorphic to $A\# \B{T}$ for every $a\in Z(A^\times)$.
\end{Theorem}

\begin{proof}
 The result follows from the fact that the image of $\varphi$ (as $\Bbbk$-vector spaces)
 is the direct sum
 \[ A \otimes \Bbbk[x] \oplus \bigoplus_{n=0}^\infty \big\langle
 a\sigma^{-1}(a)\cdots\sigma^{-n}(a)\big\rangle\otimes \operatorname{Span}_{\Bbbk}\bigl(x^{-n-1}\bigr),
 \] where $\langle
 u\rangle$ denotes the two sided ideal in $A$ generated by an element $u\in A$.
\end{proof}

In specific cases, the fact that GWAs are subalgebras of smash products was already
known~\cite[Lemma~2.3]{Boddington:ResolutionOfTypeD}. However, to the best of our knowledge,
the fact that one gets an isomorphism when the distinguished element $a\in A$ is a unit,
even though it implicitly follows from this embedding, is not fully taken advantage of in
the literature.

From now on, we identify $W_{a,\sigma}$ with $im(\varphi)$ in $A\# \B{T}$.

\subsection{Localizations of smash products with tori}
\label{sect:localizationOfSmashProducts}

Let $A$ be an algebra with a fixed automorphism $\sigma\in \operatorname{Aut}(A)$. Assume
$R\colon \B{T}\otimes A\to A\otimes \B{T}$ is the distributive law given in
equation~\eqref{monoid-dist-law}. Let $S\subseteq Z(A)$ be any multiplicative submonoid
which stable under the action of $\sigma$. The proof of the following lemma is routine
verification, and therefore, is omitted.

\begin{Lemma}\label{lem:nc-localization}
 Any multiplicative monoid $S$ in $Z(A)$ which is $\sigma$-stable is a right Ore subset
 in~$A\# \B{T}$, and $S^{-1}(A\# \B{T}) = S^{-1}A\# \B{T}$.
\end{Lemma}

\subsection{Localizations of GWAs}

As before, assume $A$ is a unital associative algebra, $a\in Z(A)$ and $\sigma\in \operatorname{Aut}(A)$.
Recall that by Theorem~\ref{thm:UnitCase} we identified the GWA $W_{a,\sigma}$ with the
subalgebra of the smash biproduct $A\# \B{T}$ generated by the algebra $A$ and the
elements $x$ and $ax^{-1}$. Then we have a tower of algebra extensions of the form
$A\# \Bbbk[x]\subset W_{a,\sigma} \subseteq A\# \B{T}$.

\begin{Theorem}\label{thm:localization}
 Consider the set $S\subset Z(A)$ of the elements of the form $\sigma^m(a^n)$, where
 $n\in\B{N}$ and~${m\in\B{Z}}$. Then the embedding of algebras
 $W_{a,\sigma}\subseteq A\# \B{T}$ is a birational equivalence with respect to the Ore
 set generated by $S$.
\end{Theorem}

\begin{proof}
 Now, by Lemma~\ref{lem:nc-localization}, we have that $S^{-1}(A\# \B{T}) = S^{-1}A\# \B{T}$,
 and by Theorem~\ref{thm:UnitCase}, we~see that the algebra $S^{-1}A\# \B{T}$ is itself
 generated by $S^{-1}A$, $x$ and $ax^{-1}$ since $a\in S^{-1}A$ is now a~unit.
\end{proof}

\subsection{Highest weight modules of GWAs}\label{sect:highest-weights}

Assume $A$ is unital associative with a distinguished element $a\in Z(A)$ and an
automorphism $\sigma\in \operatorname{Aut}(A)$. Let $V$ be a representation over the GWA $W_{a,\sigma}$.
We have an (not necessarily exhaustive) increasing filtration of submodules of the form
\[
 V^{[\ell]} = \big\{ v\in V\mid v \tl a\sigma^{-1}(a)\cdots\sigma^{-\ell}(a) = 0\big\}
\]
defined for $\ell\in\B{N}$. Let us also define
\[ V^{[\infty]} = \bigcup_{\ell\geq 0} V^{[\ell]}. \] We define $ht_{a,\sigma}(V)$ \emph{the
 height of $V$} as the smallest integer $\ell$ such that $V^{[\ell]} = V^{[\infty]}$, and
if no such integer exists we set $ht_{a,\sigma}(V)=\infty$.

Assume $V$ is a finite dimensional representation. Then $h=ht_{a,\sigma}(V)$ is necessarily
finite. Furthermore, if the height filtration satisfies $V^{[h]} = V$, then we get the
analogue of a \emph{highest weight module} for the GWA $W_{a,\sigma}$. Approaches for such
cases can be seen in~\cite{DGO:WeightModules, Meinel:DufloTheorem}.

\begin{Proposition}\label{prop:highest-weight}
 Let $S\subseteq Z(A)$ be the subset of elements of the form $\sigma^n(a^m)$ with $n\in\B{Z}$ and~${m\in\B{N}}$, and let $S^{-1}(W_{a,\sigma})$ be the localization of $W_{a,\sigma}$ at $S$. Assume
 $V$ is an arbitrary $W_{a,\sigma}$-module, and let $h=ht_{a,\sigma}(V)$. Then
 $S^{-1}V := V\otimes_{W_{a,\sigma}} S^{-1}(W_{a,\sigma})$ is isomorphic to $S^{-1}\bigl(V/V^{[h]}\bigr)$.
\end{Proposition}

\begin{proof}
 We consider the following short exact sequence of $W_{a,\sigma}$-modules:
 \[ 0 \to V^{[h]}\to V \to V/V^{[h]} \to 0 \] and use the fact that the functor
 $S^{-1}(\ \cdot\ )$ is exact.
\end{proof}

\subsection{Morphisms of algebra extensions}

An algebra $C$ together with a subalgebra $A$ is called an algebra extension. Given two
extensions~${A\subseteq C}$ and $A\subseteq C'$ of a fixed algebra $A$, a morphism $f\colon
(C,A)\to (C',A)$ of extensions is a commutative triangle of algebra morphisms of the form
\[
 \xymatrix{
 C \ar[rr]^f & & C'.\\
 & A\ar[lu]\ar[ru]
 }
\]

\subsection{Isomorphisms of smash products with tori}\label{subsect:IsoProblem}

In this subsection, we consider the isomorphism problem for smash products with
$\B{T} = \Bbbk\big[x,x^{-1}\big]$ since all isomorphism problems for GWAs birationally reduce to
isomorphism problems for such smash products.

\begin{Theorem}\label{thm:OuterAutos}
 Assume $\sigma$ and $\eta$ are two algebra automorphisms of $A$. Then the algebra
 extensions $A\subseteq A\#_{R(\sigma)} \B{T}$ and $A\subseteq A\#_{R(\eta)} \B{T}$ are
 isomorphic if and only if $\eta = u\sigma^{\pm 1}u^{-1}$ for some~${u\in A^\times}$.
\end{Theorem}

\begin{proof}
 Assume for now that $\sigma = u\eta u^{-1}$ or $\sigma = u\eta^{-1}u^{-1}$. Consider an
 arbitrary $v\in A$. In the first case, define
 $\delta\colon A\#_{R(\sigma)} \B{T}\to A\#_{R(\eta)}\B{T}$ by letting $\delta(x) = ux$
 and we get
 \begin{align*}
 \delta(x v)
 = & u x v = u \eta(v) x = \sigma(v) u x = \delta(\sigma(v)x),
 \end{align*}
 which implies $\delta$ is an isomorphism of smash biproducts. The proof for the second case is
 similar, and therefore, is omitted. On the opposite direction, assume
 $\delta\colon A\#_{R(\sigma)} \B{T}\to A\#_{R(\eta)} \B{T}$ is an isomorphism of algebra
 extensions. The one easily see that $\delta$ restricted $\B{T}$ yields an algebra monomorphism,
 and therefore, $\delta(x) = ux^{\pm 1}$ for some $u\in A^\times$ and $\delta$ restricted to $A$ is
 identity. Thus $\sigma = u\eta^{\pm 1}u^{-1}$ as expected.
\end{proof}

Notice that given an automorphism $\sigma\in \operatorname{Aut}(A)$ and its inverse $\sigma^{-1}$ extended
to $A\#_{R(\sigma)}\B{T}$ are now an inner automorphisms. From this perspective
Theorem~\ref{thm:OuterAutos} says that given two automorphism $\sigma$ and $\eta$, they
define two different smash products if their outer automorphism classes are different. In
particular, we have the following result.

\begin{Corollary}
 If $\sigma\in \operatorname{Aut}(A)$ is an inner automorphism, then the smash biproduct
 $A\#_{R(\sigma)} \B{T}$ is isomorphic to the direct product $A\times \B{T}$.
\end{Corollary}

\section{Homology of GWAs}\label{sect:Homology}

\subsection{Homology of smash products with tori}

Assume $M$ is a $\B{T}$-module $M$ via an automorphism $\sigma\colon M\to M$. We define
$
M^{\B{T}}=\{m\in M\mid t\tr m = m\}
$.
We also consider $M\otimes\B{T}$ as a $\B{T}$-bimodule via
$t\tr\bigl(m\otimes t^i\bigr) = \sigma(m)\otimes t^i$,
 $\bigl(m\otimes t^i\bigr)\tl t = m\otimes t^{i+1}
$. We denote this bimodule by $M\# \B{T}$. Now, we have the following lemma.

\begin{Lemma}\label{lem:LocalizedCoefficients}
 \begin{equation*}
 H_n(\B{T},M\#\B{T}) \cong
 \begin{cases}
 M_{\B{T}}\#\B{T} & \text{if } n=0,\\
 M^{\B{T}}\#\B{T} & \text{if } n=1,\\
 0 & \text{otherwise.}
 \end{cases}
 \end{equation*}
\end{Lemma}
\begin{proof}
 Using the fact that $\B{T}$ has Hochschild dimension 1, we can immediately conclude that $H_n(\B{T},M\#\B{T})=0$ for $n\geq 2$. As for degree 0, we have
 $H_0(\B{T},M\#\B{T}) \cong (M\#\B{T})_{\B{T}}$ and
 \[ (M\#\B{T})_{\B{T}}
 = (M\#\B{T})/[\B{T},M\# \B{T}]=\frac{M\#\B{T}}{\big\{\sigma^i(m)\otimes t^{i}-m\otimes t^i\mid m\in M,\, i\in\B{Z} \big\}}
 \cong M_{\B{T}}\# \B{T}
 \]
 As for degree 1,
 \[ H_1(\B{T},M\#\B{T})
 = \frac{\big\{m\otimes t^j\otimes t^n\mid m\otimes t^{n+j} = \sigma^{n}(m)\otimes t^{n+j}\big\}}
 {\big\{m\otimes t^{i+j}\otimes t^n - m\otimes t^i\otimes t^{n+j} + \sigma^n(m)\otimes t^{n+i}\otimes t^j \big\} }.
 \]
 Thus a copy of $M^{\B{T}}\#\B{T}$ generated by elements of the form $m\otimes t^n\otimes t$, where
 $\sigma(m)=m$ in the first homology. In the quotient, for all $m\in M$ with $\sigma^n(m)=m$, we
 have
 \begin{align*}
 m\otimes 1\otimes t^n
 = & \sigma^{n-1}(m)\otimes t^{n-1}\otimes t + m\otimes t\otimes t^{n-1},
 \end{align*}
 and by inductively repeating we get
 \begin{align*}
 m\otimes 1\otimes t^n
 = & \sum_{i=0}^{n-1} \sigma^i(m)\otimes t^{n-1}\otimes t
 \end{align*}
 for those pairs $m\otimes t^n\in M\#\B{T}$ with $\sigma^n(m)=m$. This means in the first homology
 terms of the form $m\otimes t^i\otimes t^n$ with $\sigma^n(m)=m$ are replaced by terms of the form
 \[ \Biggl(\sum_{j=0}^{n-1}\sigma^j(m)\Biggr)\otimes t^{n+i-1}\otimes t .
 \]
 Notice that the sum on the left is $\sigma$-invariant. Result follows.
\end{proof}

Thanks to Lemma~\ref{lem:LocalizedCoefficients}, we now have the following result by using
 the spectral sequence given in Proposition~\ref{base2}.
\begin{Proposition}\label{prop:torus-action}
 Let $\sigma\in \operatorname{Aut}(A)$ and assume $\sigma$ acts on $\CH_*(A)$ diagonally extending the
 action on $A$. Let $\CH_*(A)_{\B{T}}$ and $\CH_*(A)^{\B{T}}$ respectively be the complex
 of coinvariants and invariants of~$\sigma$. Then
 \[ HH_n(A\# \B{T})
 \cong H_n(\CH_*(A)_{\B{T}})\otimes \B{T} \oplus H_{n-1}\bigl(\CH_*(A)^{\B{T}}\bigr)\otimes \B{T}.
 \]
\end{Proposition}

\subsection{Algebraic and separable endomorphisms}

We call an algebra endomorphism $\sigma\in \operatorname{End}(A)$ \emph{algebraic} if there is a
polynomial $f(t)\in\Bbbk[t]$ such that $f(\sigma) = 0$ in $\operatorname{End}(A)$. For an algebraic
endomorphism $\sigma$ of $A$, the monic polynomial $f(t)$ with the minimal degree that
satisfies $f(\sigma)=0$ is called \emph{the minimal polynomial of $\sigma$.} We call an
algebraic endomorphism $\sigma\in \operatorname{End}(A)$ as \emph{separable} if the minimal polynomial of
$\sigma$ is separable.

Notice that all endomorphisms of a finite dimensional $\Bbbk$-algebra are algebraic.
Regardless of the dimension, all automorphisms of finite order and all nilpotent non-unital
endomorphisms are also algebraic. If $\Bbbk$ has characteristic~0, automorphisms of finite
order are separable, but nilpotent non-unital endomorphisms are not.

\subsection{Algebras with separable automorphisms}\label{sect:separable}

For a fixed algebraic automorphism $\sigma\in \operatorname{Aut}(A)$, let $\operatorname{Spec}(\sigma)$ be the set of
unique eigenvalues of~$\sigma$, and let $A^{(\lambda)}$ be the $\lambda$-eigenspace of
$\sigma$ corresponding to $\lambda \in \operatorname{Spec}(\sigma)$.

\begin{Theorem}\label{thm:separable}
 Assume $\sigma\in \operatorname{Aut}(A)$ is separable with minimal polynomial $f(x)$, and let
 $B$ be the quotient $\Bbbk[x]/\langle f(x)\rangle$. Then
 \[ H_n(A\# \B{T})
 = H_n\bigl(\CH_*^{(1)}(A)\bigr)\otimes\B{T} \oplus H_{n-1}\bigl(\CH_*^{(1)}(A)\bigr)\otimes \B{T}
 \] and
$ H_n(A\# B) = H_n\bigl(\CH_*^{(1)}(A)\bigr)\otimes B$, where $\CH_*^{(1)}(A)$ is generated by homogeneous tensors of the form
$a_0\otimes\cdots\otimes a_n$ with $ a_i\in A^{(\lambda_i)}
$ and $\lambda_1\cdots\lambda_n=1
$ for every $n\geq 0$.
\end{Theorem}

\begin{proof}
 One can extend the distributive law $R\colon \B{T}\otimes A\to A\otimes\B{T}$ given in
 equation~\eqref{monoid-dist-law} to a~distributive law of the form
 $R\colon B\otimes A\to A\otimes B$. Notice that since $f(x)$ is separable, $B$ is a~product of a~finite number of copies of $\Bbbk$, and therefore, is $0$-smooth. Then the result for $A\# B$
 immediately follows from Corollary~\ref{base3}. On the other hand, the result for $A\# \B{T}$
 follows from Proposition~\ref{prop:torus-action}.
\end{proof}

Note that Theorem~\ref{thm:separable} solves the smoothness problem for smash products with
$\B{T}$, and therefore \emph{the birational smoothness} problem for all GWAs, provided that
the action is implemented via a separable automorphism. Namely, a smash product with
$\B{Z}$ via a separable automorphism is smooth if and only if the complex subcomplex of
invariants $\CH_*^{(1)}(A)$ has bounded homology. In the next subsection, we solve the
birational smoothness problem for all GWAs without requiring automorphism to be separable.

\subsection{Localization of GWAs in homology}\label{sect:localizationGWA}

Consider the set $S$ of elements of the form $\sigma^m(a^n)$ in $Z(A)$ where $n\in\B{N}$ and
$m\in\B{Z}$. Let $\Bbbk[ S]$ be the (commutative) subalgebra of $A$ generated by $S$,
and let $S^{-1}\Bbbk[ S]$ be its localization at~$S$. Then we have that
$S^{-1}A=A\otimes_{\Bbbk[ S]}S^{-1}\Bbbk[ S]$. Now let
$\Bbbk[ S]_{\B{T}}$ be the algebra of coinvariants of~$\Bbbk[ S]$ which is
given by the following quotient:
\[ \Bbbk[ S]_{\B{T}} := \frac{\Bbbk[ S]}{\langle \sigma(s)-s\mid
 s\in S\rangle}. \]

\begin{Corollary} \label{cor:localization}
 We have
 \begin{align*}
 HH_n\bigl(S^{-1}W_{a,\sigma}\bigr)
 \cong{} & HH_n\bigl(S^{-1}A\# \B{T}\bigr)\\
 \cong{} & H_n\bigl(\CH_*(A)_{\B{T}}\otimes_{\Bbbk[ S]_{\B{T}}}S^{-1}(\Bbbk[ S]_{\B{T}})\bigr) \otimes \B{T}
 \oplus H_{n-1}\bigl(\CH_*\bigl(S^{-1}A\bigr)^{\B{T}}\bigr)\otimes \B{T},
 \end{align*}
 where we view $\CH_*(A)$ as an $\Bbbk[ S]$-module and $\Bbbk[ S]_{\B{T}}$-module on the
 coefficient.
\end{Corollary}

\begin{proof}
 By Theorem~\ref{thm:localization}, we have $S^{-1}W_{a,\sigma} \cong S^{-1}A\# \B{T}$. Now, we
 consider the algebra extension~$S^{-1}A\subseteq \bigl(S^{-1}A\bigr)\# T$ for which by \cite{KaygSutl20}
 there is a spectral sequence whose first page is
 \[ E^1_{p,q} = H_q\bigl(S^{-1}A,\CH_p\bigl(S^{-1}A\# \B{T}|S^{-1}A\bigr)\bigr) = H_q\bigl(S^{-1}A,\CH_p\bigl(\B{T},S^{-1}A\# \B{T}\bigr)\bigr) \]
 that converges to $HH_*\bigl(S^{-1}A\# \B{T}\bigr)$. Since $S\subseteq Z(A)$, by
 \cite{Brylinski:Localizations}, we know that
 \[ E^1_{p,q}\cong H_q\bigl(A,S^{-1}\CH_p\bigl(\B{T}, S^{-1}A\# \B{T}\bigr)\bigr)\cong H_q\bigl(A,\CH_p\bigl(\B{T},S^{-1}A\# \B{T}\bigr)\bigr). \] Thus we have an isomorphism of the form $HH_*\bigl(S^{-1}A\# \B{T}\bigr)\cong
 H_*\bigl(A\# \B{T},S^{-1}A\# \B{T}\bigr)$. Then by Lemma~\ref{lem:LocalizedCoefficients} and Proposition~\ref{prop:torus-action}, we get
 \[ HH_n\bigl(S^{-1}W_{a,\sigma}\bigr)\cong H_n\bigl(\CH_*\bigl(A,S^{-1}A\bigr)_{\B{T}}\bigr)\otimes\B{T} \oplus
 H_{n-1}\bigl(\CH_*\bigl(S^{-1}A\bigr)^{\B{T}}\bigr)\otimes\B{T}.
 \] Since $S\subseteq Z(A)$, we get that $\CH_*\bigl(A,S^{-1}A\bigr) = S^{-1}\CH_*(A)$. On the other hand,
 both the coinvariants functor $(\ \cdot\ )_{\B{T}}$ and localization functor $S^{-1}(\ \cdot\ )$ are
 specific colimits, and colimits commute. Then
 \[ \CH_n\bigl(S^{-1}W_{a,\sigma}\bigr)\cong S^{-1}\CH_n(A)_{\B{T}}\otimes\B{T} \oplus \CH_*\bigl(S^{-1}A\bigr)^{\B{T}}\otimes\B{T}.
 \]
 On the other hand,
 \begin{gather*}
 \bigl(S^{-1}\CH_*(A)\bigr)_{\B{T}} \cong \bigl(\CH_*(A)\otimes_{\Bbbk[S]} S^{-1}\Bbbk[S]\bigr)_{\B{T}} \cong \CH_*(A)\otimes_{\Bbbk[S]\rtimes\B{T}} S^{-1}\Bbbk[S] \\ S^{-1}\CH_*(A)_{\B{T}} \cong \CH_*(A)_{\B{T}}\otimes_{\Bbbk[S]_{\B{T}}} S^{-1}(\Bbbk[S]_{\B{T}})
 \end{gather*}
 since the $S$-localization of an $A$-module $M$ can be written as
 \[ S^{-1}M = M\otimes_A S^{-1}A \cong M\otimes_{k[S]} S^{-1}\Bbbk[S].\tag*{\qed}
 \]\renewcommand{\qed}{}
\end{proof}

\section{Homology calculations}\label{sect:Calculations}

\subsection{The rank-1 Weyl algebra}\label{sect:WeylAlgebra}

Consider the ordinary Weyl algebra as a GWA as we did in Section~\ref{sect:OrdinaryWeyl}. Now,
let $S$ be the multiplicative system generated by elements of the form $(t-m)$ where $m\in\B{Z}$.
Since $t$ is not a~unit in $A$ we see that $W_{t,\sigma}$ is the proper subalgebra of
$\Bbbk[t]\# \B{T}$ generated by $x$ and $tx^{-1}$ where the distributive law is defined as
$R(x\otimes t) = (t-1)\otimes x$. Since there is no non-constant rational function invariant under
the action $\sigma(f(t))=f(t-1)$, we get that \smash{$\CH_*\bigl(S^{-1}\Bbbk[t]\bigr)^{\B{T}} = \CH_*(\Bbbk)$}. Next,
we see that the subalgebra generated by $S$ is $A=\Bbbk[t]$ itself. Moreover, since $\sigma(t)-t=1$
we get that $\Bbbk[S]_{\B{T}}$ is zero, and therefore, we get
\[
 HH_n\bigl(S^{-1}A_1\bigr) =
 \begin{cases}
 \B{T} & \text{if } n=1,\\
 0 & \text{otherwise}
 \end{cases}
\]
for every $n\geq 0$.

\subsection[The enveloping algebra U(G sl\_2)]{The enveloping algebra $\boldsymbol{U(\G{sl}_2)}$}\label{sect:Usl2}

The universal enveloping algebra of $\G{sl}_2$ is given by the presentation
\[ \frac{\Bbbk\{E,F,H\}}{\langle EH-(H-2)E, FH - (H+2)F, EF-FE-H \rangle}. \] The center
of this algebra is generated by the Casimir element
\[ \Omega = 4 FE + H(H + 2) = 4EF + H(H - 2). \] In this subsection, we would like to write
a generalized Weyl algebra isomorphic to $U(\G{sl}_2)$.

Let $A=\Bbbk[c,t]$ and $a=c-t(t+1)$. Define $\sigma$ to be the algebra automorphism
defined by~${\sigma(f(c,t)) = f(c,t-1)}$ for every $f(c,t)\in A$. In this case,
$W_{a,\sigma}$ is generated by $c$, $t$, $x$ and~${(c-t(t+1))x^{-1}}$ in the smash product
algebra $A\# \B{T}$. The GWA $W_{a,\sigma}$ is isomorphic to~$U(\G{sl}_2)$ via an
isomorphism defined as
$H\mapsto 2t$, $ E\mapsto x$, $ F\mapsto (c-t(t+1))x^{-1}$,
see~\cite[Example~2.2]{FarinatiSolotarAlvarez:GeneralizedWeyAlgebras}.

Let us define $S$ to be the multiplicative system generated by elements of the form
$c-(t-n)(t-n-1)$, for $ n\in\B{Z}$.
Then $S^{-1}W_{t,\sigma}$ is isomorphic to $S^{-1}\Bbbk[c,t]\# \B{T}$, and
\smash{$\CH_*\bigl(S^{-1}A\bigr)^{\B{T}}=\CH_*(\Bbbk[c])$}. Moreover, the subalgebra of $A=\Bbbk[c,t]$ generated by
$S$ is $A$ itself and since $\sigma(t)-t = 1$, we again get that~${\Bbbk[S]_{\B{T}}=0}$.
Therefore,
\[
 HH_n\bigl(S^{-1}U(\G{sl}_2)\bigr) =
 \begin{cases}
 \Bbbk[c]\otimes \B{T} & \text{if } n=1,2,\\
 0 & \text{otherwise}.
 \end{cases}
\]

\subsection[Primitive quotients of U(G sl\_2)]{Primitive quotients of $\boldsymbol{U(\G{sl}_2)}$}

One can also consider $B_\lambda := W_{a,\sigma}/\langle c-\lambda\rangle$, where
$W_{a,\sigma}$ is $U(\G{sl}_2)$ as we defined above. These algebras are also GWAs since we
can realize them using $A=\Bbbk[t]$, $a=\lambda-t(t+1)$ with the $\sigma$ given by
$t\mapsto t-1$. See~\cite[Section 3]{BavulaJordan:IsoProblem}.

In this case, using a similar automorphism we used for $U(\G{sl}_2)$, we can replace $S$ with the
multiplicative system generated by elements of the form $\mu-(t-n)$ and $\mu+(t-n)$, where~${\mu\in\Bbbk}$ is fixed and $n$ ranges over $\B{Z}$. Then $\Bbbk[S] = \Bbbk[t]$ and
$S^{-1}B_\lambda\cong S^{-1}\Bbbk[t]\# \B{T}$. In this case, \smash{$\CH_*\bigl(S^{-1}\Bbbk[t]\bigr)^{\B{T}}$} is~$\CH_*(\Bbbk)$ and $\Bbbk[S]_{\B{T}}=0$ since $\sigma(t)-t=1$ as before. Then we get
\[ HH_n\bigl(S^{-1}B_\lambda\bigr)\cong
 \begin{cases}
 \B{T} & \text{if } n=1,2,\\
 0 & \text{otherwise}
 \end{cases}
\] for every $n\geq 0$.

\subsection{Quantum 2-torus}\label{subsect:quantumTorus}

Fix an element $q\in\Bbbk^\times$ which is not a root of unity. Let $A =\Bbbk\big[t,t^{-1}\big]$ and
let $a=t$ as in the case of the ordinary Weyl algebra. But this time, let us define
$\sigma\in \operatorname{Aut}(A)$ to be the algebra automorphism given by $\sigma(f(t)) = f(qt)$ for every
$f(t)\in A$. The smash biproduct algebra $A\# \B{T}$ is the algebraic quantum 2-torus
$\B{T}^2_q$ and the GWA $W_{a,\sigma}$ is the quantum torus itself since $a=t$ is a unit.

Note that for every $u\in A$ and $m\in\B{Z}$ we have $\sigma^m(u)\neq u$ unless $m=0$
since $q$ is not a root of unity. Thus \smash{$\CH_*(A)_\B{T} = \CH_*(A)^\B{T} = \CH^{(0)}_*(A)$}
where
\[
 \CH^{(0)}_m(A)
 = \operatorname{Span}_{\Bbbk}\biggl(t^{n_0}\otimes\cdots\otimes t^{n_m}\mid n_1,\ldots,n_m\in\B{Z}
\text{ with } 0 = \sum_i n_i \biggr),
\]
which gives us just the group homology of $\B{Z}$. Then by
Proposition~\ref{prop:torus-action}, we get
\[ HH_n\bigl(\B{T}^2_q\bigr) \cong \Bbbk^{\binom{2}{n}}\otimes \B{T}
\] for every $n\geq 0$ as expected.

\subsection[The quantum enveloping algebra U\_q(Gsl\_2)]{The quantum enveloping algebra $\boldsymbol{U_q(\G{sl}_2)}$}\label{sect:QuantumSL2}

For a fixed $q\in\Bbbk^\times$, the quantum enveloping algebra of the lie algebra $\G{sl}_2$
is given by the presentation
\[ \frac{\Bbbk\big\{K,K^{-1},E,F\big\}}{\big\langle KE - q^2EK, KF - q^{-2}FK,
 EF-FE=\frac{K-K^{-1}}{q-q^{-1}}\big\rangle}.
\] As before, we assume $q$ is not a root of unity. There is an element $\Omega$ in the
center of $U_q(\G{sl}_2)$ called the \emph{quantum Casimir element} defined as
\[
 \Omega = EF + \frac{q^{-1}K + qK^{-1}}{\bigl(q-q^{-1}\bigr)^2}
 = FE + \frac{qK+q^{-1}K^{-1}}{\bigl(q-q^{-1}\bigr)^2}.
\]
See~\cite[Section~I.3]{BrownGoodearl:QuantumGroups}. Our first objective is to give a GWA that
is isomorphic to $U_q(\G{sl}_2)$.

We start by setting $A=\Bbbk\big[c,t,t^{-1}\big]$ together with
$
 a = c - \bigl(q^{-1}t+qt^{-1}\bigr)
$
and $\sigma\in \operatorname{Aut}(A)$ given by $\sigma(f(c,t)) = f\bigl(c,q^2 t\bigr)$ for every
$f(c,t)\in\Bbbk\big[c,t,t^{-1}\big]$. Define an algebra map
$\gamma\colon W_{a,\sigma}\to U_q(\G{sl}_2)$ given on the generators by
\[
 t\mapsto K, \qquad
 c\mapsto \bigl(q-q^{-1}\bigr)^2\Omega, \qquad
 x\mapsto \bigl(q-q^{-1}\bigr)F,\qquad
 ax^{-1} \mapsto \bigl(q-q^{-1}\bigr)E.
\] Notice that the inverse of $\gamma$ is defined easily as
\[ K\mapsto t, \qquad
 E\mapsto \frac{ax^{-1}}{q-q^{-1}},\qquad
 F\mapsto \frac{x}{q-q^{-1}}.
\] One can show that both $\gamma$ and its inverse are well-defined by showing the
relations are preserved.

Now, let $S$ be the multiplicative system in $A$ generated by the elements of the form
\[
 c - \bigl(q^{-2n+1}t+q^{2n-1}t^{-1}\bigr)\qquad \text{for}\quad n\in\B{Z}.
\]
In this case too, the subalgebra of $A$ generated by $S$ is $A$ itself. Then we have
\[
 \CH_*\bigl(S^{-1}A\bigr)^{\B{T}} \cong \CH_*^{(0)}\bigl(\Bbbk\big[c,t,t^{-1}\big]\bigr) \cong \CH_*(A)_{\B{T}}.
\]
On the other hand, since $\sigma(t)-t = \bigl(q^2-1\bigr)t$ and $t$ is a unit, we get that
$\Bbbk\langle S\rangle_{\B{T}} = 0$. Thus, as in the case of $U(\G{sl}_2)$ we get
\begin{align*}
 HH_n\bigl(S^{-1}W_{a,\sigma}\bigr)
 \cong & HH_n\bigl(S^{-1}U_q(\G{sl}_2)\bigr)
 \cong
 \begin{cases}
 \Bbbk[c]\otimes\B{T} & \text{if } n=1,2,\\
 0 & \text{otherwise }
 \end{cases}
\end{align*}
for every $n\geq 0$.

\subsection[The quantum matrix algebra CO\_q(M\_2)]{The quantum matrix algebra $\boldsymbol{\C{O}_q(M_2)}$}

For a fixed $q\in\Bbbk^\times$ the algebra $\C{O}_q(M_2)$ of quantum $2\times 2$ matrices
is given by the presentation
\begin{gather*}
 bc = cb,\qquad
 ab = q^{-1}ba,\qquad
 ac = q^{-1}ca,\qquad
 db = qbd,\\
 dc = qcd,\qquad
 ad-da = \bigl(q^{-1}-q\bigr)bc.
\end{gather*}
The quantum determinant
$\Omega = ad - q^{-1}bc = da - qbc$ generates the center of this algebra,
see~\cite[pp.~4--8]{BrownGoodearl:QuantumGroups}.

Now, let $A = \Bbbk[u,v,w]$ with the distinguished element $u + q vw\in A$ where we set
$\sigma(f(u,v,w))\allowbreak=f\bigl(u,q^{-1}v,q^{-1}w\bigr)$ for every $f(u,v,w)\in A$. Then the GWA
$W_{a,\sigma}$ is the subalgebra of $A\# \B{T}$ generated by $A$, $x$ and
$(u + qvw)x^{-1}$, and it is isomorphic to $\C{O}_q(M_2)$ via
\[ u\mapsto\Omega,\qquad
 v\mapsto b,\qquad
 w\mapsto c,\qquad
 x\mapsto a,\qquad
 (u+qvw) x^{-1}\mapsto d,
\]
and its inverse is
\[ a\mapsto x,\qquad
 b\mapsto v,\qquad
 c\mapsto w,\qquad
 d\mapsto (u+qvw)x^{-1}.
\]

Since $\C{O}_q({\rm GL}_2)$ is obtained by localizing $\C{O}_q(M_2)$ at the quantum determinant, we see
that $\C{O}_q({\rm GL}_2)$ is isomorphic to $u^{-1}W_{a,\sigma}$ which itself is a GWA with $A$ replaced
by $\Bbbk\big[u,u^{-1},v,w\big]$ with the remaining datum unchanged.

On the other hand, $\C{O}_q({\rm SL}_2)$ is the quotient of $\C{O}_q(M_2)$ by the two sided
ideal generated by~$u-1$, and therefore, is again a GWA with the same datum where this
time we replace $A$ by~$\Bbbk[u,v,w]/\langle u-1\rangle$. We also know that
$\C{O}_q({\rm GL}_2)$ is isomorphic (as algebras only) to $\C{O}_q({\rm SL}_2)\times \Bbbk[\Omega]$.

For the remaining of the section, we are going to concentrate on $\C{O}_q({\rm SL}_2)$ only given
as the subalgebra of $\Bbbk[v,w]\# \B{T}$ generated by $v$, $w$, $x$ and
$(1+qvw)x^{-1}$.

Now, let $S$ be the Ore set generated by elements of the form $1+q^{2n+1}vw$ for $n\in\B{Z}$. Then~$S^{-1}\C{O}_q({\rm SL}_2)$ is isomorphic to $S^{-1}\Bbbk[v,w]\# \B{T}$. In this case, since $q$ is not a
root of unity, we get that
\smash{$\CH_*\bigl(A,S^{-1}A\bigr)^{\B{T}}=\CH_*(\Bbbk) = \CH_*(A)_{\B{T}}$}. The subalgebra of $\Bbbk[v,w]$
generated by $S$ is the polynomial algebra $\Bbbk[vw]$ over the indeterminate~$vw$. Since
$\sigma(vw)-vw=\bigl(q^{-2}-1\bigr)vw$, we get that $\Bbbk[vw]_{\B{T}} = \Bbbk$. Hence
\[
 HH_n\bigl(S^{-1}\C{O}_q({\rm SL}_2)\bigr) \cong \Bbbk^{\binom{2}{n}}\otimes\B{T}
\] for every $n\geq 0$.

\subsection[Quantum group CO\_q(SU\_2)]{Quantum group $\boldsymbol{\C{O}_q({\rm SU}_2)}$}\label{sect:OqSU2}

Let us fix $q\in\Bbbk^\times$. The algebraic quantum group $\C{O}_q({\rm SU}_2)$ is the
noncommutative $*$-algebra generated by two non-commuting indeterminates $s$ and $x$ subject
to the following relations:
\[
 x^*x = 1 - s^*s, \qquad x x^* = 1 - q^2 s^*s, \qquad s^*s = s s^*,\qquad xs =
 qsx,\qquad xs^*=qs^*x.
\]
See \cite[p.~4]{Hadfield:Podles}. One can write $\C{O}_q({\rm SU}_2)$ as a GWA $W_{a,\sigma}$ by
letting $A=\Bbbk[s,s^*]$ with the distinguished element $a\in A$ is defined as $1-s^*s$ and
$\sigma(f(s,s^*))= f(qs,qs^*)$ for every $f(s,s^*)\in\Bbbk[s,s^*]$.

Let $S$ be the multiplicative system in $A$ generated by elements of the form $q^{2n}s^*s-1$ for~${n\in\B{Z}}$. Then $S^{-1}\C{O}_q({\rm SU}_2)$ is isomorphic to $S^{-1}A\# \B{T}$ by
Theorem~\ref{thm:localization}. If we assume that $q\in \Bbbk^\times$ is not a~root of unity, we get
that
\smash{$
 \CH_*\bigl(S^{-1}A\bigr)^{\B{T}} = \CH_*(\Bbbk) = \CH_*(A)_{\B{T}}
$}
We also see that the subalgebra of $\Bbbk[s,s^*]$ generated by $S$ is the polynomial algebra
$\Bbbk[ss^*]$, and since $\sigma(ss^*)-ss^* = \bigl(q^2-1\bigr)ss^*$ we get that
$\Bbbk\langle S\rangle_{\B{T}} = \Bbbk$. Then
\[
 HH_n\bigl(S^{-1}\C{O}_q({\rm SU}_2)\bigr) \cong \Bbbk^{\binom{2}{n}}\otimes\B{T}
\]
for every $n\geq 0$.

\subsection{Podle\'s spheres}\label{sect:PodlesSphere}

For a fixed $q\in\Bbbk^\times$, the algebra of functions $\C{O}_q\bigl(S^2\bigr)$ on standard Podle\'s
quantum spheres~\cite{Hadfield:Podles,Podl87} is the subalgebra of $\C{O}_q({\rm SU}_2)$
generated by elements $s^*s$, $xs$ and $s^*x^*$. This means $\C{O}_q\bigl(S^2\bigr)$ is the
subalgebra of the smash product $\Bbbk[s,s^*]\# \B{T}$ generated by the elements $s^*s$,
$sx$ and $s^*(1-s^*s)x^{-1}$. One can give a presentation for the Podle\'s sphere as
\[
xt = q^2tx, \qquad yt = q^{-2} ty, \qquad yx = -t(t-1), \qquad xy =
-q^2t\bigl(q^2t-1\bigr),
\]
 then we get a GWA structure if we let $A = \Bbbk[t]$ and where we set
$t=s^*s$ with $a=-t(t-1)$ and $\sigma(f(t)) = f\bigl(q^2 t\bigr)$ for every $f(t)\in A$.

Let $S$ be the multiplicative system in $A$ generated by the set $\big\{ t\bigl(t-q^{2n}\bigr)\mid n\in\B{Z} \big\}$,
then ${S^{-1}\C{O}_q\bigl(S^2\bigr) \cong S^{-1}A\# \B{T}}$. Instead of this generating set, one can use
$\{t\}\cup \big\{\bigl(t-q^{2n}\bigr)\mid n\in\B{Z}\big\}$ to get the same localization. Then we get that $\Bbbk[S]$
is $A$ itself. If we assume that $q\in \Bbbk^\times$ is not a~root of unity, we get that
\smash{$
 \CH_*\bigl(S^{-1}A\bigr)^{\B{T}} = \CH^{(0)}_*\bigl(\Bbbk\big[t,t^{-1}\big]\bigr)$}, and $ \CH_*(A)_{\B{T}} = \CH_*(\Bbbk)$.
In this case, $\Bbbk[S]_{\B{T}} = \Bbbk$ since $\sigma(t)-t = \bigl(q^2-1\bigr)t$. Thus
\[ HH_n\bigl(S^{-1}\C{O}_q\bigl(S^2\bigr)\bigr)
 \cong \Bbbk^{\binom{2}{n}}\otimes \B{T}
\]
for every $n\geq 0$.

\subsection{Parametric Podle\'s spheres}

In~\cite{Hadfield:Podles}, Hadfield defines another family of Podle\'s spheres
$\C{O}_{q,c}\bigl(S^2\bigr)$ given by a presentation is equivalent to the following:
\[
 x t = q^2 t x, \qquad
 x^* t = q^{-2} t x^*,\qquad
 x^*x = c - t(t - 1),\qquad
 x x^* = c - q^2 t\bigl(q^2 t -1 \bigr).
\]
If we set $A = \Bbbk[c,t]$, and let the distinguish element $a\in A$ be $c-t(t-1)$
together with $\sigma(f(c,t)) = f\bigl(c,q^2t\bigr)$ for every $f(c,t)\in A$, we get a GWA structure
on $\C{O}_{q,c}\bigl(S^2\bigr)$ similar to the GWA structure on $U(\G{sl}_2)$ where we changed only
the algebra automorphism from $\sigma(f(c,t)) = f(c,t-1)$ to $\sigma(f(c,t)) = f\bigl(c,q^2t\bigr)$.

Let $S$ be the multiplicative system in $A$ generated by the elements of the form
$c-q^{2n}t\bigl(q^{2n}t-1\bigr)$. If we assume that $q\in \Bbbk^\times$ is not a root of unity, we
conclude that
\[
 \CH_*\bigl(S^{-1}A\bigr)^{\B{T}} = \CH_*(\Bbbk[c]) = \CH_*(A)_{\B{T}}
\]
which allows us to conclude
\[
 HH_n\bigl(S^{-1}\C{O}_{q,c}\bigl(S^2\bigr)\bigr)
 \cong \Bbbk^{\binom{2}{n}}\otimes\Bbbk[c]\otimes\B{T}
\]
for every $n\geq 0$.

\subsection*{Acknowledgments}

We would like to thank the anonymous referees for their careful reading of the paper, and
suggestions and corrections they provided as it greatly improved the main arguments and
presentation of the paper. This work was completed while the author was on academic leave
at Queen's University from Istanbul Technical University. The author is supported by the
Scientific and Technological Research Council of Turkey (T\"UB\.{I}TAK) sabbatical grant 2219.
The author would like to thank both universities and T\"UB\.{I}TAK for their support.

\pdfbookmark[1]{References}{ref}
\LastPageEnding

\end{document}